\def\draft{y}

\documentclass{amsart}
\usepackage{amssymb, amsmath, latexsym, amscd,graphicx}
\usepackage{fullpage}

\def\printname#1{
	\if\draft y
		\smash{\makebox[0pt]{\hspace{-0.5in}
			\raisebox{8pt}{\tt\tiny #1}}}
	\fi }

\newtheorem{theorem}{Theorem}[section]

\newtheorem{lemma}[theorem]{Lemma}
\newtheorem{claim}[theorem]{Claim}
\newtheorem{proposition}[theorem]{Proposition}

\newtheorem{fact}[theorem]{Fact}

\theoremstyle{definition}
\newtheorem{definition}[theorem]{Definition}

\newcommand{\A}{\mathcal A}

\newcommand{\bl}{\mathcal B^e_n}

\title{Homotopy approximations to the space of knots, Feynman diagrams,
and a conjecture of Scannell and Sinha}
\author{James Conant}
\email{jconant@math.utk.edu}
\address{Department of Mathematics, University of Tennessee, Knoxville, TN, 37996}

\begin{document}
\begin{abstract}Scannell and Sinha considered a spectral sequence to calculate the
rational homotopy groups of spaces of long knots in $\mathbb R^n$, for $n\geq 4$. 
At the end of the paper they conjecture that when $n$ is odd, the terms on the antidiagonal at the $E^2$ stage precisely give the space of irreducible Feynman diagrams related to the theory of Vassiliev invariants. 
In this paper we prove that conjecture. This has the application that the path components of the terms of the Taylor tower for the space of long knots in $\mathbb R^3$ are in one-to-one correspondence with quotients of the module of Feynman diagrams, even though the Taylor tower does not actually converge. This provides strong evidence that the stages of the Taylor tower give rise to universal Vassiliev knot invariants in each degree. Our proof yields a sequence of new presentations for the space of irreducible Feynman diagrams.
\end{abstract}
\maketitle
\large
\section{Introduction}
Consider the space of long knots, $\operatorname{Emb}(\mathbb R,\mathbb R^n)$, which are a fixed line outside of a compact set. 
According to the calculus of Goodwillie and collaborators, one can define  homotopy-theoretic approximations to the space of knots:
$$ev_k\colon \operatorname{Emb}(\mathbb R,\mathbb R^n)\to AM_k.$$
The map $ev_k$ induces isomorphisms on homology and homotopy to a larger and larger extent as $k$ increases, provided that $n\geq 4$. In the classical case of $n=3$, we still get knot invariants:
$$\pi_0(ev_k)\colon \pi_0(\operatorname{Emb}(\mathbb R,\mathbb R^3))\to\pi_0(AM_k).$$
We conjectured in \cite{bcss} that these maps are actually universal Vassiliev invariants of degree $k-1$ over the integers. The calculations in this paper will show that 
$\pi_0(AM_k)$ is a quotient of the space of primitive Feynman diagrams of degree $k$ that appear in the theory of Vassiliev invariants of knots. (It is a quotient because in the spectral sequence calculations, higher differentials might possibly kill off some of the space.)

To establish this conjecture, we analyze Scannell and Sinha's spectral sequence computations in 
\cite{ss}. They consider a spectral sequence which converges to the rational homotopy groups of the space of long knots, constructed via the Taylor approximations $AM_k$, or more precisely, via equivalent cosimplicial models.
The main result of this paper is that, when $n$ is odd, the terms along the antidiagonal of the $E^2$ page are isomorphic to spaces $\A^I_k$ of primitive Feynman diagrams.

These spaces $\A^I_k$ are known to rationally classify primitive Vassiliev invariants of degree $k$ up to lower-degree invariants. The bulk of the current paper is devoted to
giving an alternate presentation for $\A^I_k$, which is hopefully of independent interest. The usual presentation is via trivalent graphs attached to a line segment, modulo the $\operatorname{STU}, \operatorname{IHX},\operatorname{AS}$ and $\operatorname{SEP}$ relations. 
The new presentation is via connected trivalent trees attached to a line segment modulo the so-called $\operatorname{STU}^2$, $\operatorname{IHX}$ and $\operatorname{AS}$ relations. The $\operatorname{STU}^2$ relation is pictured in Figure~\ref{blash}.

In the last section we show that the $E^1$ terms on the antidiagonal are isomorphic to the space of trees attached to a directed line segment, modulo $\operatorname{IHX}$ relations. On the antidiagonal, passing to the $E^2$ page involves dividing by the image of the differential, which has the effect of introducing $\operatorname{STU}^2$ relations, thus completing the argument.

The main theorem of this paper (Theorem \ref{main}) was proven independently and with a completely different approach by Lambrechts and Tourtchine \cite{lt}.  

{\bf Acknowledgments:} This paper arose from discussions with Dev Sinha, and forms part of a joint project with him, Ryan Budney and Kevin Scannell. I'd like to thank them for their helpful discussions. I'd also like to thank Ted Stanford for helpful discussions about alternate presentations of the space of primitive Feynman diagrams, Viktor Tourtchine for inspiration and helpful discussion, and the referee for pointing out many typographical errors in an earlier version. This research was partially supported by NSF grant DMS 0305012.

\section{An overview of Scannell and Sinha's result}\label{overview}
Let $\bl$ be the free graded Lie algebra generated by
elements
$x_{ij}$ of degree $1$ where $1\leq i,j\leq n$, subject to the following relations:
\begin{align*}
&x_{ij}=-x_{ji}\\
&x_{ii}=0\\
&[x_{ij},x_{lm}]=0 \text{ if }\{i,j\}\cap\{l,m\}=\emptyset\\
&[x_{ij},x_{jl}]=[x_{jl},x_{li}]=[x_{li},x_{ij}]
\end{align*}

Let $M_{d,n}$ be the submodule of the degree $d$ summand of $\bl$ generated by brackets of
elements $x_{in}$ where all $i<n$ appears as an index. This definition is equivalent to saying that $M_{d,n}$ is the submodule of the degree $d$ summand of $\bl$ generated by brackets of elements $x_{ij}$ where all of the indices $1,\ldots, n$ appear. (This follows from \cite[Algorithm 5.2]{ss}, which will convert an iterated bracket of generators $x_{ij}$ to a sum of brackets involving only generators $x_{in}$ and brackets which don't involve the index $n$. These latter terms cannot arise if all indices are present.)

There is a differential $d\colon M_{d,n}\to M_{d,n+1}$ given by
$$d=\sum_{i=0}^{n+1}(-1)^i\partial^i$$ where

$\partial^l(x_{ij})= x_{\sigma^l(i)\sigma^l(j)}$ if $i,j\neq l$ and
$\partial^l(x_{ij})= x_{i\sigma^l(j)}+x_{i+1 \sigma^l(j)}$ if $i=l$, where
$\sigma^l(i) = i$ if $i<l$ and equals $i+1$ if $i>l$.

\begin{theorem}[Scannell and Sinha]
Let $k\geq 4$ be even.
There is a spectral sequence which converges, over $\mathbb Q$, to 
$$\pi_*(\operatorname{Emb}(I,\mathbb R^k\times I))$$
whose $E^1$ term is given by $E^1_{-n,d(k-1)+1}=M_{d,n}$ and whose $d^1$ is given by the differential $d$ defined above.
\end{theorem}

In fact, a recent result \cite{ALTV} indicates that the spectral sequence collapses at the $E^2$ term. See Volic's survey paper \cite{volic}.

For the case of classical knots, $k=2$, there is no convergence result, but Sinha and Scannell conjectured that
the submodule of classes along the anti-diagonal correspond to primitive Vassiliev knot invariants. 



The main theorem of the current paper is the following. We will define $\A^I_n$ in the next section.
\begin{theorem}\label{main}
Let $k$ be even. Then 
$$E^2_{-(n+1),(n+1)}=M_{n,n+1}/\operatorname{im}(d)\cong \A^I_n\otimes\mathbb Q.$$
\end{theorem}

Here is a different perspective on $d$. In calculating $\partial^l$ one takes each instance of $l$ and replaces it by either an $l$ or an $l+1$. Thus there are $2^k$ terms in $\partial^l(c)$ where $k$ is the number of times the index $l$ appears in $c$.

\begin{proposition}
The differential $d$ is equal to $\sum_{i=1}^{n} (-1)^i\tilde{\partial}^i$,
where the operator $\tilde{\partial}^l(c)$ consists of those terms in $\partial^l(c)$ in which all indices appear. 
\end{proposition}
\begin{proof}
This is straightforward.
\end{proof}

\section{Other presentations of $\A^I_n$}
\begin{figure}
\begin{center}
$\underset{(a)}{\includegraphics[width=3in]{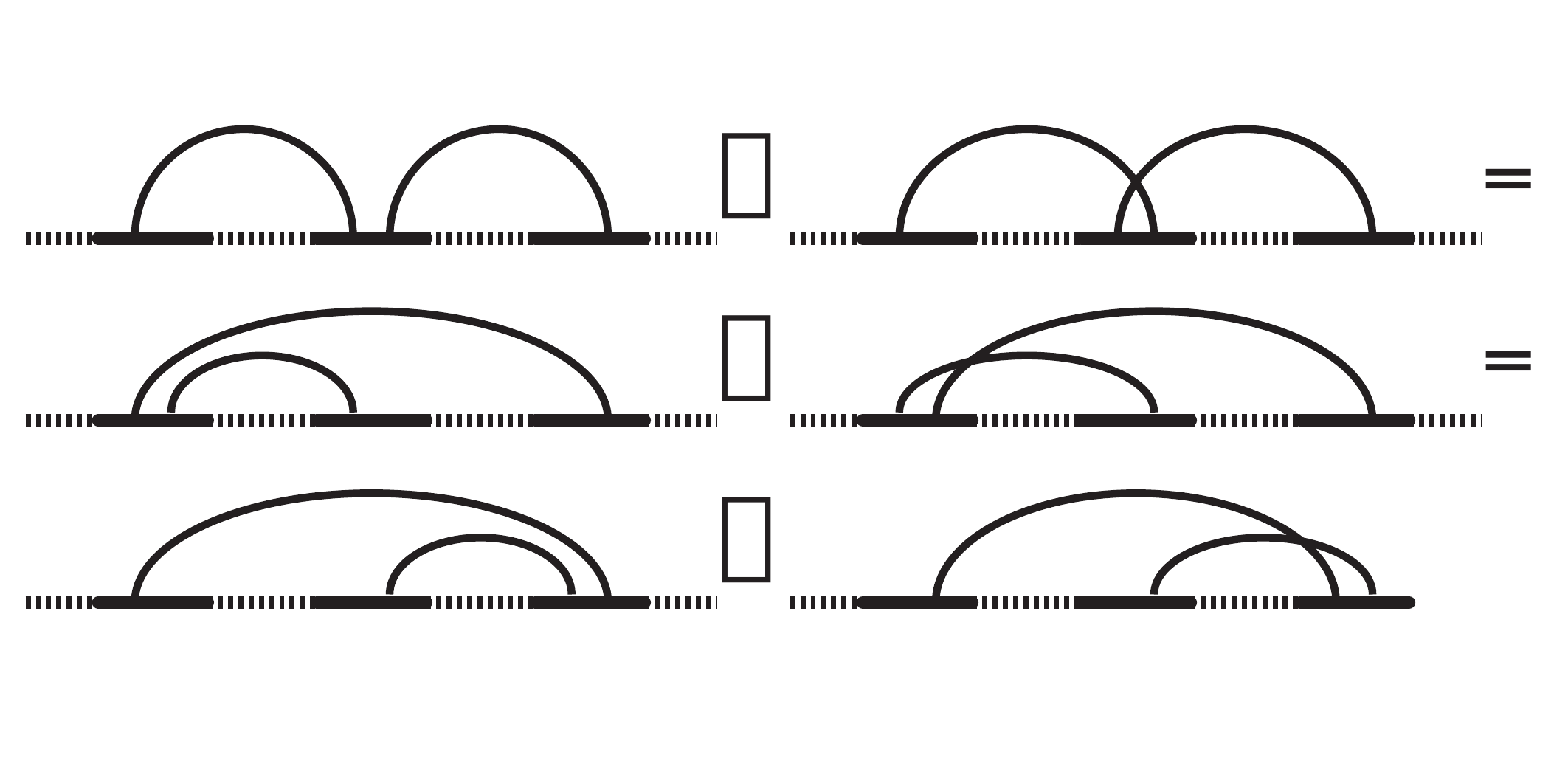}}\hfill\underset{(b)}{\includegraphics[width=2.5in]{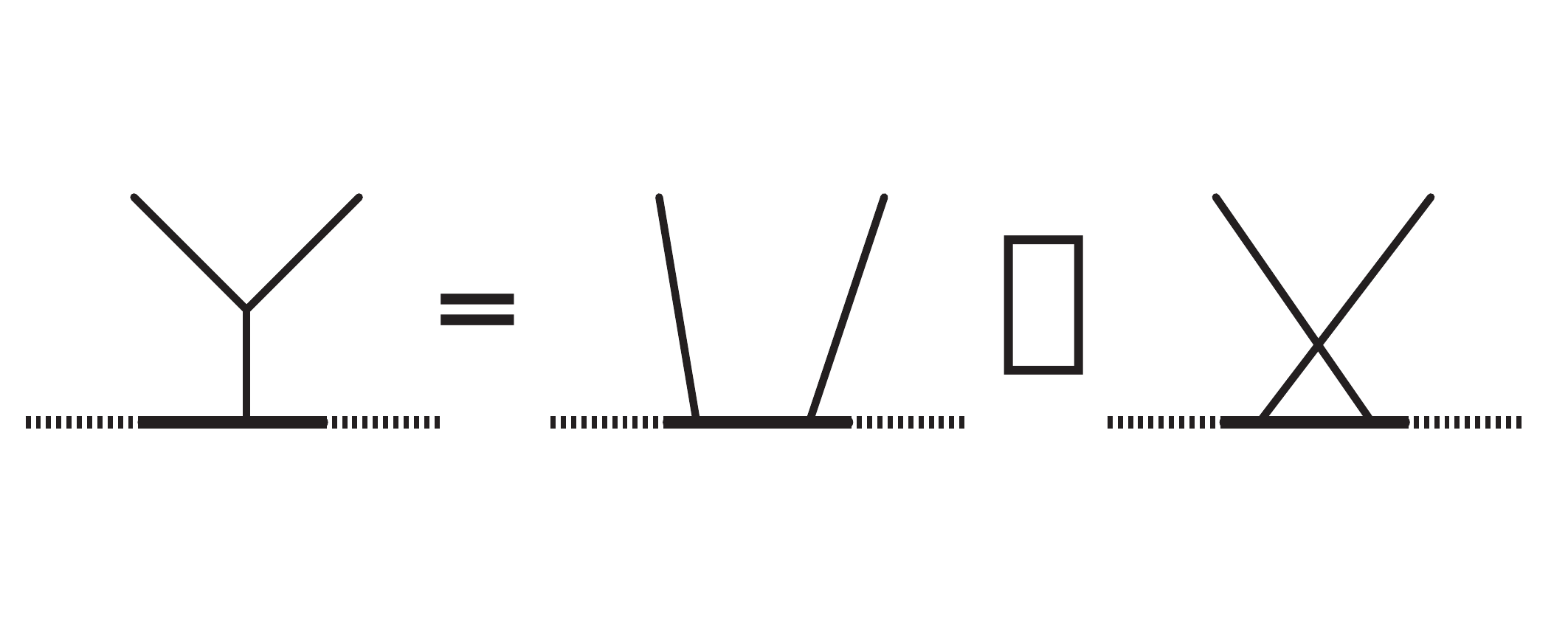}}$\\
$\underset{(c)}{\includegraphics[width=3.5in]{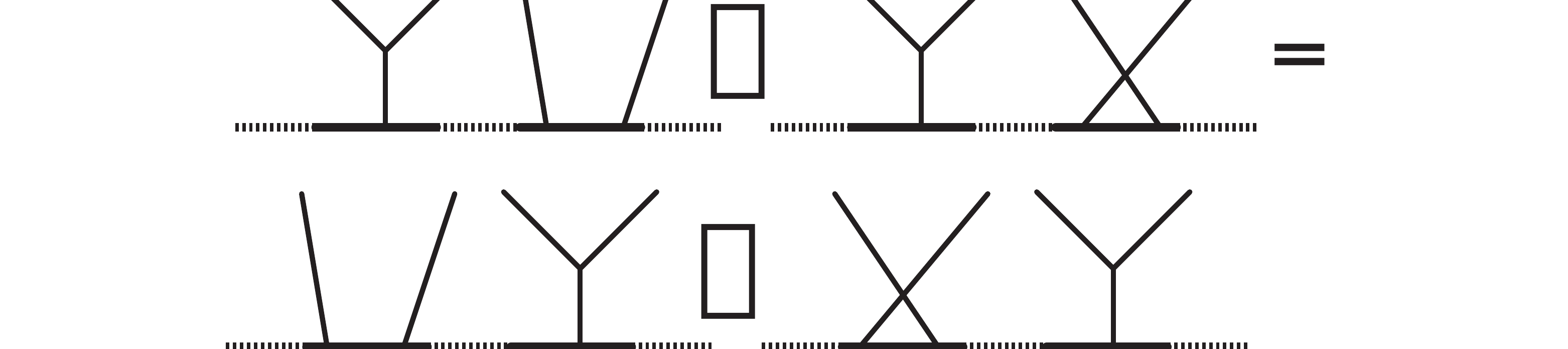}}\hfill\underset{(d)}{\includegraphics[width=3in]{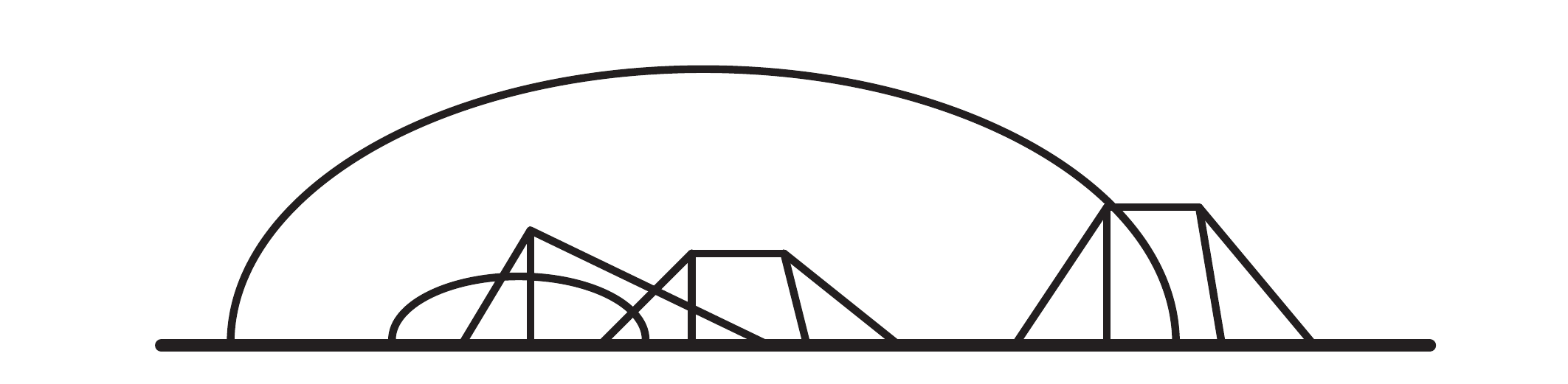}}$
\end{center}
\caption{(a) The 4T relation. (b) The STU relation. (c) The $\operatorname{STU}^2$ relation. (d) A separated diagram.} \label{blash}
\end{figure}

Here is a quick and dirty review of some spaces of diagrams related to Vassiliev invariants. 
Let $\A^I_n$ be the $\mathbb Z$-module of chord diagrams. It is generated by diagrams formed by attaching $n$ chords to a directed line segment along distinct pairs of points.
The relations are of two forms. The first is called the $4T$ relation, and is pictured in Figure~\ref{blash}. The second is called the $\operatorname{SEP}$ relation and is the relation setting separated diagrams to zero. A separated diagram (see Figure \ref{blash}d) is one which there is an isolated clump of trees which do not interact with the rest of the trees. More formally, thinking of the diagram as being immersed in the plane, there is a circle which intersects the line segment in two points, does not intersect any of the trees, and contains trees in both its exterior and interior. 
The superscript $I$ in $\mathcal A^I_n$ indicates we are dividing out separated diagrams and stands for ``irreducible."

Another presentation of this same module is by attaching unitrivalent trees (or more generally graphs) to a line segment instead of just chords. The total number of vertices, including both the vertices internal to the trees and the vertices occuring where the trees attach to the line segment, is $2n$. Also each vertex internal to a tree has a specified cyclic order. The relations are of three forms. The so-called AS relation says that switching the cyclic order at a vertex is the same as multiplication by $-1$. The STU relation is as in Figure~\ref{blash}b. The SEP relation sets separated diagrams to $0$ as before.
In this context, elements of $\A^I_n$ are often called Feynman diagrams (on a directed line segment).

Let's begin by defining some modules closely related to $\A^I_n$.

\begin{definition}

\noindent\begin{enumerate}
\item
$\tilde{A}^I_{k,n}$ is defined to be the module of degree $n$, vertex-oriented, Feynman diagrams on a directed line segment, such that deletion of the line segment yields $k$ tree components.
\item $\A^I_{k,n}$ is the quotient of $\tilde{A}^I_{k,n}$ by the relations
\begin{enumerate}
\item $\operatorname{IHX}$
\item $\operatorname{STU^2}$
\item $\operatorname{SEP}$.
\end{enumerate}
Here $\operatorname{IHX}$ is the standard relation which takes place away from the directed line segment. $\operatorname{STU^2}$ refers to the relation below: one takes an element (called a \emph{template}) of $\tilde{A}^I_{k-1,n}$ or a template which has $k-1$ tree components and a component which is a unitrivalent graph homotopy equivalent to a circle, and breaks a vertex open via the standard $\operatorname{STU}$ relation in two different ways. (In this latter case, one needs to break open vertices that will convert the circle-like graph into a tree.)
  Pictured in Figure \ref{blash} is a relation which breaks two vertices apart, but one could also break the same vertex apart using two outgoing edges. Finally $\operatorname{SEP}$ sets any separated diagram equal to $0$.
\end{enumerate}
\end{definition}

Notice that $\A^I_{n,n}$ is the usual module of chord diagrams, proven by Bar-Natan\cite{bn} to coincide with general Feynman diagrams modulo STU and IHX. So $\A^I_{n,n}\cong\A^I_n$. Note that $\operatorname{STU^2}$ coincides with the $\operatorname{4T}$ relation in this case.

We will use the following fact repeatedly: 

\begin{fact} If $k\leq n-1$, then the $\operatorname{STU}^2$ relations are generated by $\operatorname{STU}^2$ relations which involve breaking apart two distinct vertices of its template. 
\end{fact}
\begin{proof}
In the case $k\leq n-1$ the template has at least $2$ trivalent vertices. Any $\operatorname{STU}^2$ relation involving the same vertex twice can be rewritten as the difference of two $\operatorname{STU}^2$ relations, each involving two distinct vertices.
\end{proof}

Now we turn to the main theorem of this section:

\begin{theorem}
$\A^I_{k,n}\cong \A^I_n$ for all $1\leq k\leq n$ with the possible exception of $k=n-1$.
\end{theorem}

\begin{proof}
Here is the strategy of the proof. We will define maps 
$$\Psi_k\colon\A^I_{k,n}\to \A^I_{k-1,n}\text{ and  }\Phi_{k-2}\colon\A^I_{k-2,n}\to \A^I_{k-1,n}$$ which will make sense for $k\neq n$. The fact that $\Phi_{k-1}\circ\Psi_k=\operatorname{Id}$ will then imply that each $\Psi_k$, $k<n-1$ is injective. 
We will also show that each $\Psi_k$ is onto. At this point we can then conclude that $\A^I_{n-2,n}\cong \A^I_{n-3,n}\cong\cdots\cong \A^I_{1,n}$. To bridge the last gap, we observe that even though $\Psi_n$ is not well-defined, it is well-defined on the module $\tilde{\A}^I_n$, where no relations are present:
 $\tilde{\Psi}_n\colon \tilde\A^I_{n,n}\to \A^I_{n-1,n}$. We then argue that the composition
 $\Psi_{n-1}\circ\tilde\Psi_n$ does in fact kill the submodule of relations, giving rise to a well-defined map
 $\Psi^\prime\colon \A^I_{n,n}\to\A^I_{n-2,n}$, which is evidently onto. In a similar vein we define
 $\Phi^\prime=\tilde\Phi_{n-2}\circ\Phi_{n-1}$, and argue that $\Phi^\prime\circ\Psi^\prime=\operatorname{Id}$, implying that $\Psi^\prime$ is injective, completing the proof.

The map $\Phi_k$ is easier to define. Choose any trivalent vertex which is connected by an edge to the directed line segment, and apply an STU relation to get a difference of two elements in $\A^I_{k+1,n}$. 

\begin{claim} $\Phi_k$ is well-defined, $k\neq n-2$. 
\end{claim}
\begin{proof}
If one expands a different vertex, then the result is related by an $\operatorname{STU^2}$-relation.  The fact that $\Phi(\operatorname{SEP})\subset \operatorname{SEP}$ is obvious. As for $\Phi(\operatorname{STU^2})$,  we may assume that there are two distinct vertices involved, as pictured in Figure~\ref{blash}. In each of these four terms, calculate $\Phi$ by applying an STU relation to
the visible trivalent vertex.
 The result is zero on the nose. 
 The fact that $\Phi(\operatorname{IHX})\subset \operatorname{IHX}$ is where we use the hypothesis $k\neq n-2$. If $k=n-1$, there are no IHX relations present, and the result trivially holds. If $k<n-2$, then there are at least three trivalent vertices, and we can split apart one not involved in the IHX relation when calculating $\Phi$, the result obviously lying in the IHX subspace.
\end{proof}


Now we define maps $\Psi_k\colon \A^I_{k,n}\to\A^I_{k-1,n}$, $k\neq n$.
They are defined in the following way. Given a diagram, $C$, in $\tilde\A^I_{k,n}$, let $C_1$ denote the component whose leg hits the directed line segment first. (Farthest to the left.) $-\Psi(C)$ is a sum of diagrams where each leg of $C_1$ is consecutively attached to the legs of the rest of the diagram, moving from right to left until the leg is planted to the left of the rest of the diagram, whereupon the process is repeated for the next leg. See Figure~\ref{psidef}. One way to think of the oval notation is that
one chooses a horizontal slice of the grey region, transverse to the trees inside of it, and distributes the endpoint of the edge to all of the edges hitting this slice. It is not hard to show that this is well-defined modulo $\operatorname{IHX}$ relations. However, an easy way to get a well-definition is to take the horizontal slice to be close to the directed line segment.

Observe that $\Psi$ is defined so that it becomes the identity once one introduces $\operatorname{STU}$ relations, as it represents the difference between having $C_1$ attached in its original position and having $C_1$ slid all the way to the left (which is zero modulo $\operatorname{SEP}$.)

\begin{figure}
\begin{center}
\includegraphics[width=4in]{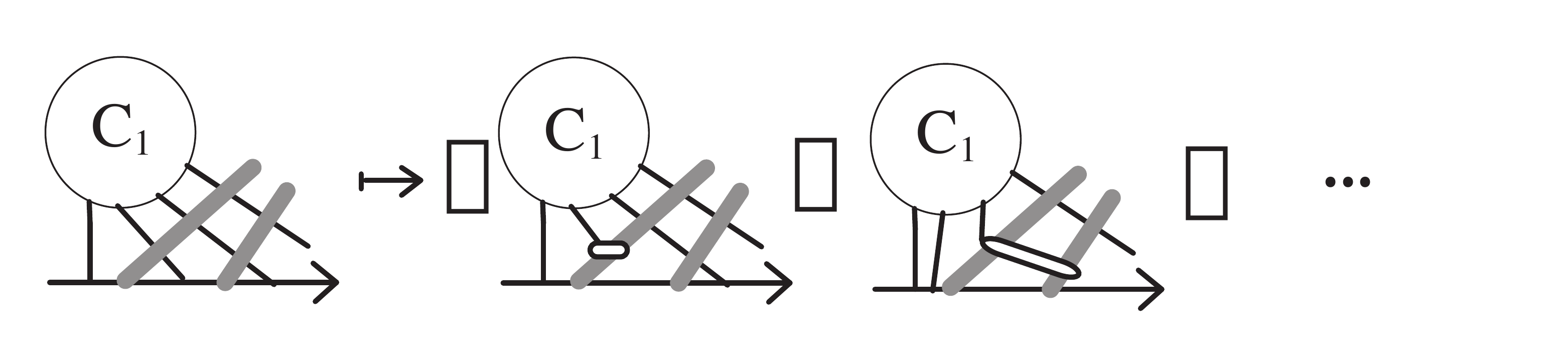}
\caption{The definition of $\Psi$. The circle represents the trivalent tree $C_1$ which is attached first to the line segment. The thick grey lines represent regions where legs of other trees may attach. A leg attaching to grey regions by a white oval means a sum of diagrams where the leg attaches to each of the legs inside the grey region.\label{psidef}}
\end{center}
\end{figure}

\begin{claim}\label{clm}
The following equation holds on the diagrammatic level. That is the domain of $\Psi$ is not divided by any relations, although the range \emph{is}.\\

\centerline{\noindent\includegraphics[width=3in]{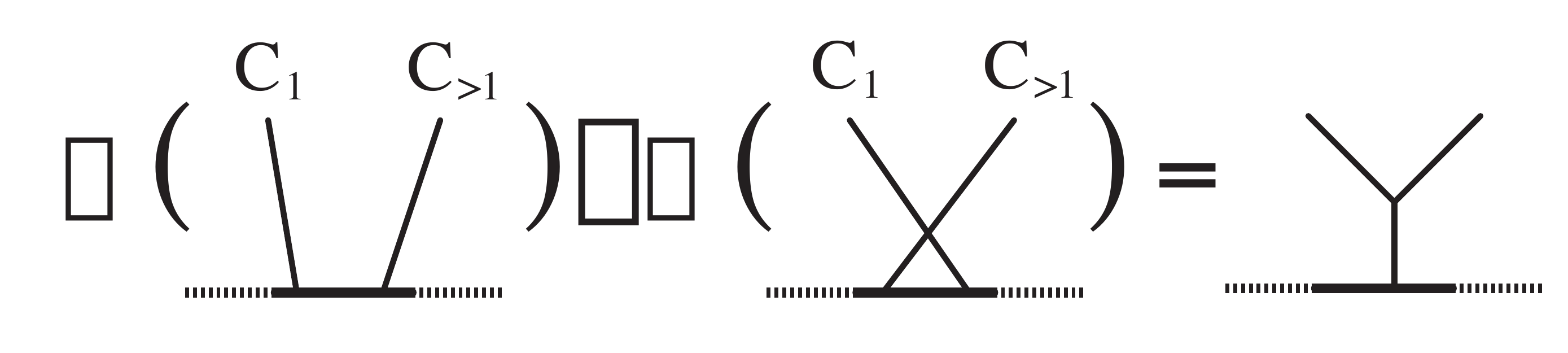}}

Here the leg labeled $C_1$ is part of the tree $C_1$, whereas the other leg is part of a different tree.
\end{claim}
\begin{proof}
Let $L$ be the visible leg of $C_1$. Call the legs of $C_1$ attached to the left of $L$ \emph{prior legs}.
Applying $\Psi$ to the second term, eventually we get a term formed by crossing the leg $L$ across the leg to its left, creating a term similar to the one on the right, except that all prior legs of $C_1$ have been shifted to the left. After this, all the terms from each of the left-hand terms above actually match, and so cancel in pairs. Up to this point, the terms
  differ by a transposition of two legs. Using an $\operatorname{STU^2}$ relation, turn this transposition into a trivalent vertex in exchange for splitting apart the vertex created by $\Psi$. This gives a sum of terms which are easily seen to combine to give the difference between having the prior legs shifted to the left and having them in their original position.
\end{proof}

\begin{claim}
$\Psi$ is well-defined.
\end{claim}
\begin{proof}
That is, it vanishes on $\operatorname{SEP}+\operatorname{IHX}+\operatorname{STU^2}$. For $\operatorname{SEP}$, consider an isolated clump of trees not containing $C_1$. As we drag the legs of $C_1$ across this clump of trees, the result does not obviously lie in the $\operatorname{SEP}$ subspace. However, note that the sum of attaching a leg to the legs of an isolated clump of trees is $0$ modulo $\operatorname{IHX}$. (See the proof of Lemma 3.1 in \cite{bn}.) Thus the non-$\operatorname{SEP}$ terms cancel out.

 It is also routine to see that $\Psi(\operatorname{IHX})\subset \operatorname{IHX}$. Finally we consider $\Psi(\operatorname{STU^2})$. There are two cases:

\vspace{1em}

\noindent{\bf Case 1.} The $\operatorname{STU^2}$ relation does not involve the first component.
 Without loss of generality, assume our $\operatorname{STU^2}$ involves two distinct vertices.
All of the summands of $\Psi$ which don't involve a leg interacting with a leg in the two active sites, denoted $\alpha$ and $\beta$,
of the $\operatorname{STU^2}$ relation can be grouped to lie in $\operatorname{STU^2}$. Thus we need only consider dragging a leg across one of the active sites, say $\alpha$. The top of Figure~\ref{t1t6} depicts the summands of $\Psi$ where a leg, $L$, of $C_1$ interacts with the legs in the $\alpha$ region.

\begin{figure}[ht!]
\begin{center}
\includegraphics[width=.7\linewidth]{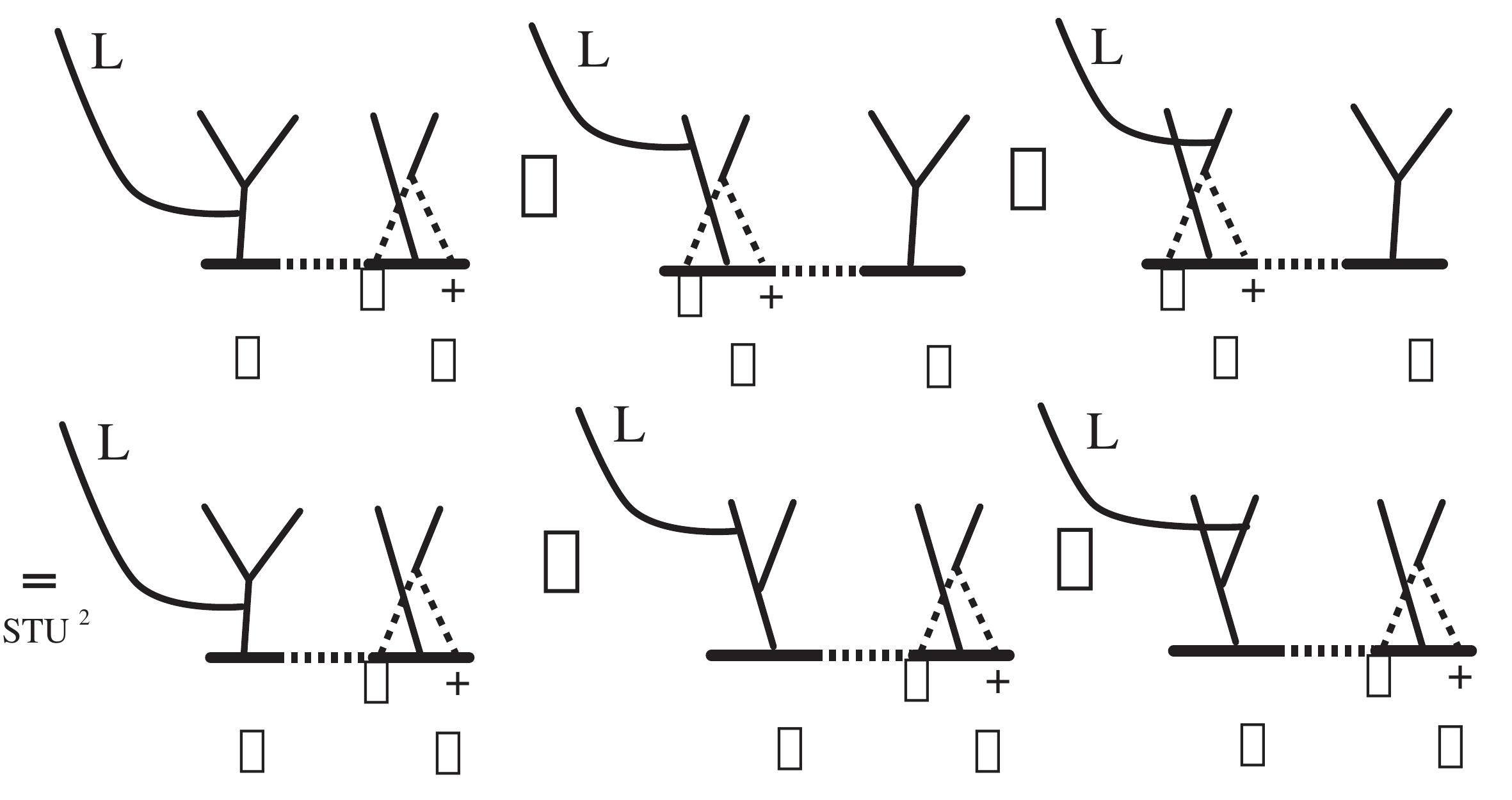}
\caption{From the proof of Case 1. The dotted lines represent a sum of two terms in which the solid line attaches to the bottom line segment in the two indicated ways, with the sign depicted.}\label{t1t6}
\end{center}
\end{figure}

Using $\operatorname{STU}^2$ relations on the last two terms, we get a difference of two IHX relators, which is therefore trivial, as shown in the bottom of Figure~\ref{t1t6}.

\vspace{1em}

 \noindent{\bf Case 2.} The $\operatorname{STU^2}$ relation does involve the first component.
Here there are three subcases. 

{\bf Subcase 1.} The $\operatorname{STU^2}$ relation comes from a template graph with a loop being split apart in two ways. Thus all of the feet involved in the $\operatorname{STU}^2$ relation come from the component $C_1$. Thus one can picture the rest of the graph being slid to the right of $C_1$ as opposed to $C_1$ being slid to the left, and then apply the argument from Case 1. 

{\bf Subcase 2.}  The relation comes from splitting apart two distinct vertices of a template tree, $T$. In this case, the component $C_1$ changes. Applying $\operatorname{STU}$ to \emph{both} vertices we get a union of three trees, $T_1$, $T_2$ and $T_3$. The $\operatorname{STU}^2$ relation comes in two pairs of terms, where in each pair two of the trees, $T_i$, are spliced together and a foot of this resulting tree differs by a transposition with the foot of the third tree in the two terms of the pair.
Then Claim~\ref{clm} indicates that $\Psi$ has the effect of converting each pair of terms to the original template tree $T$, with opposite signs, so that the total is zero.

{\bf Subcase 3.} The relation comes from splitting apart two distinct vertices, each on separate trees, $T_1$ and $T_2$. Assume that $T_1$ has a foot attaching farthest to the left.
Using Claim~\ref{clm}, $\Psi$ applied to two of the terms in which $T_1$ is split apart
just gives the template $T$. In the other two terms, we are sliding $T_1$ 
to the left, creating new trivalent vertices along the way, to calculate $\Psi$.
Suppose first that the site where $T_2$ is being split is to the right of all of the legs of
$T_1$. Then a similar argument to Claim~\ref{clm} allows us to destroy the created trivalent vertices in exchange for reassembling $T_2$. The sum of terms as a result is the difference between attaching $T_1$ all the way to the left (which is zero modulo $\operatorname{SEP}$) and attaching $T_1$ in the original position (which is the template $T$).
Now we consider the case when the splitting site for $T_2$ is to the left of some legs of $T_1$. For the above argument to work, we need to show that the sum of terms where a leg of $T_1$ attaches to the two legs in $T_2$'s splitting site is equivalent modulo $\operatorname{STU}^2$ to the difference of two terms where $T_2$ is reassembled, and the leg of $T_1$ attaches just before and just after $T_2$'s splitting site, as pictured in Figure~\ref{t4}.
\begin{figure}[ht!]
\begin{center}
\includegraphics[width=.7\linewidth]{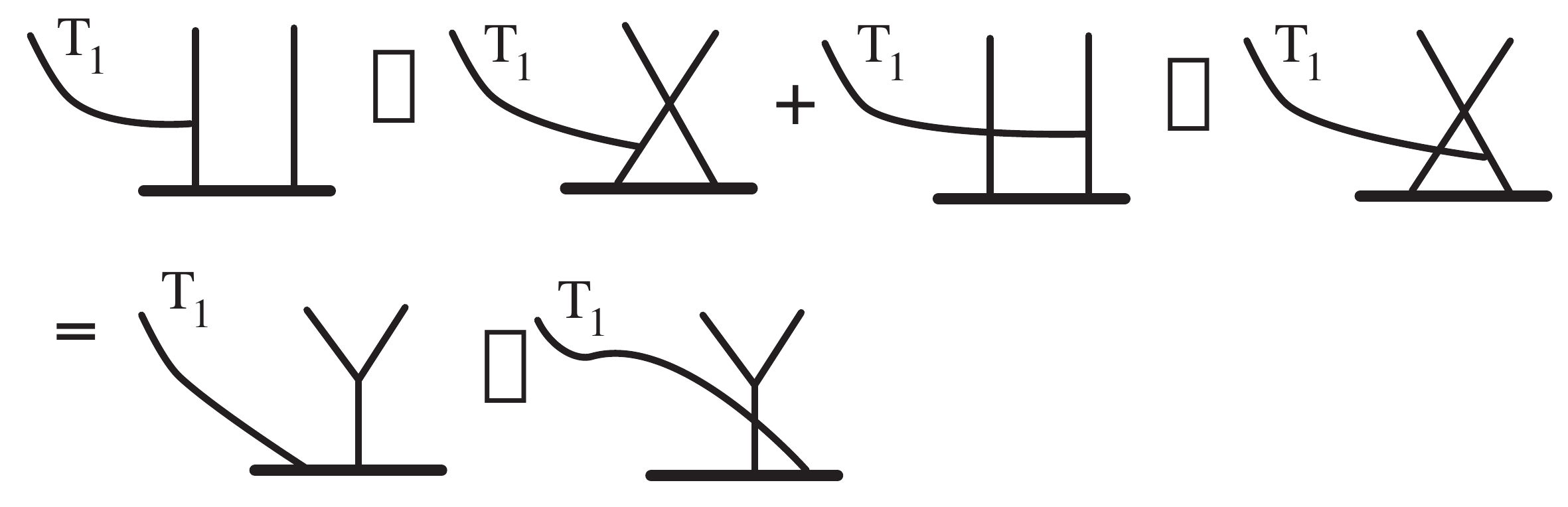}
\caption{From the proof of Subcase 3.}\label{t4}
\end{center}
\end{figure}
To see that the equality in Figure~\ref{t4} holds, use $\operatorname{STU}^2$ relations to combine the diagrams in pairs, creating new trivalent vertices, in exchange for splitting apart some other trivalent vertex, such as the one on $T_1$ that we know exists. Equality of the two sides is now an expression of the $\operatorname{IHX}$ relation.
\end{proof}

\begin{claim}
$\Phi\circ\Psi = \operatorname{Id}$. Hence $\Psi$ is injective.
\end{claim}
\begin{proof}
This is fairly obvious. $\Psi$ is calculated by dragging feet of $C_1$ to the left picking up trivalent vertices as the feet attach to the intervening legs.
To calculate $\Phi$ pull apart these very trivalent vertices. The result is the difference between the original diagram and the one gotten by sliding $C_1$ all the way to the left, and this latter term lies in $\operatorname{SEP}$.\end{proof}

Let $$\Phi^\prime\colon \A^I_{n-2,n}\to \A^I_{n,n}$$ be given by pulling apart two vertices adjacent to the line segment via $\operatorname{STU}$.
Let $$\Psi^\prime\colon \A^I_{n,n}\to \A^I_{n-2,n}$$ be defined by $\Psi\circ\tilde{\Psi}$ where $\tilde{\Psi}\colon \tilde{\A}^I_{n,n}\to\A^I_{n-1,n}$ is the map $\Psi$ defined on the diagrammatic level.
The fact that $\Phi^\prime$ is well-defined is straightforward.
\begin{claim}
$\Psi^\prime$ is well-defined. 
\end{claim}
\begin{proof}
We need to verify that 4T and SEP are in the kernel of $\Psi'$. Indeed SEP is already in the kernel of $\tilde{\Psi}$.
To show that 4T is in the kernel, suppose first that the first chord $C_1$ is not involved in the relation. Then all terms except the ones where $C_1$ attaches to chords in the relation are clearly in the $4T\subset STU^2$ subspace. Thus we need to show that the sum of terms where $C_1$ attach to the 4T chords are in the kernel of $\tilde{\Psi}$. Using Claim~\ref{clm}, this is evident. 
An \emph{ad hoc} calculation takes care of the case when $C_1$ is part of the $4T$ relation; this calculation is included in the appendix.
\end{proof}
\begin{claim}
$\Psi^\prime$ is injective. 
\end{claim}
\begin{proof}
This is a consequence of $\Phi^\prime\circ\Psi^\prime=\operatorname{Id}$, which follows from the following commutative diagrams:
$$
\begin{CD}
\A^I_{n,n}@>{\Psi'}>>\A^I_{n-2,n}\\
@VV{\cong}V @VVV\\
\A^I_n@>{=}>>\A^I_n\\
\end{CD}
\phantom{abcdedfg}
\begin{CD}
\A^I_{n,n}@<{\Phi'}<<\A^I_{n-2,n}\\
@VV{\cong}V @VVV\\
\A^I_n@<{=}<<\A^I_n\\
\end{CD}
$$
where the vertical arrows are induced by dividing by $\operatorname{STU}$ relations. The left-hand isomorphism in each diagram is proven to be so in \cite{bn}. Note that $\Phi'$ and $\Psi'$ were defined so that they become the identity once one divides by $\operatorname{STU}$ relations, explaining the commutativity of these diagrams. 
\end{proof}

\begin{claim}
The maps $\Psi$ and $\Psi^\prime$ are surjective.
\end{claim}
\begin{proof}
 If the first component has a trivalent vertex then Claim~\ref{clm} indicates that it is in the image of $\Psi$. 
 On the other hand, these types of diagrams actually generate $\A^I_{k,n}$ for $k<n$, which can be seen as follows. 
Given a diagram not of this form, it is easy to verify that it is equal to a sum of diagrams, where the right hand foot of the first chord is attached to each of the legs between it and the left-hand foot, and where in each of these terms a trivalent vertex has been split into a sum of two terms. 

The fact that $\Psi^\prime$ is onto follows since both $\Psi$ and $\tilde\Psi$ are onto. (Recall that $\tilde\Psi$ has a larger domain than $\Psi$, but that the range is the same.)
\end{proof}

\end{proof}

\section{Relation to $E^2$}
Now let us to return to the question of why 
$M_{n,n+1}/im(d)\cong \A^I_{1,n}\otimes\mathbb Q$.
Indeed we will now define a map
$$\Delta\colon M_{n,n+1}\to \A^I_{1,n}\otimes\mathbb Q.$$

We can think of $\A^I_{1,n}\otimes\mathbb Q$ as a vector space spanned by vertex-oriented trivalent trees with leaves that have a specific bijection with $1,\ldots, n+1$ , which represent the order in which the leaves attach to the directed line segment. The $\operatorname{IHX}$, $\operatorname{AS}$ and $\operatorname{STU^2}$ relations can then be interpreted in this context.  

We begin with some definitions
\begin{definition}
Given two sets of indices $\alpha$ and $\beta$, which have a single index, $i$, in common,
define $\alpha |\beta$ to be the number of pairs of indices $(a,b)$, of $\alpha\setminus\{i\}$ and $\beta\setminus\{i\}$ respectively where $a>b$.
\end{definition}

\begin{definition}
Let $\mathfrak A_n$ be the vector space spanned by vertex-oriented trivalent trees which have a specified injective map from the leaves to the set $\{1,\ldots, n+1\}$. That is, we can think of $\mathfrak A_n$ as spanned by trees with leaves labeled by the numbers $1,\ldots,n+1$ with no repetition. We divide by the usual AS and IHX relations.

Define a Lie bracket on $\mathfrak A_n$ by the following rule. Suppose $\alpha$ is the set of numbers labeling a tree $t_1\in\mathfrak A_n$, and $\beta$ is the set of number labeling a tree $t_2\in\mathfrak A_n$, then $[t_1,t_2]$ is defined to be zero unless $\alpha\cap\beta$ consists of a single number, say $i$, in which case it is defined by the equation:

\centerline{\includegraphics[width=.7\linewidth]{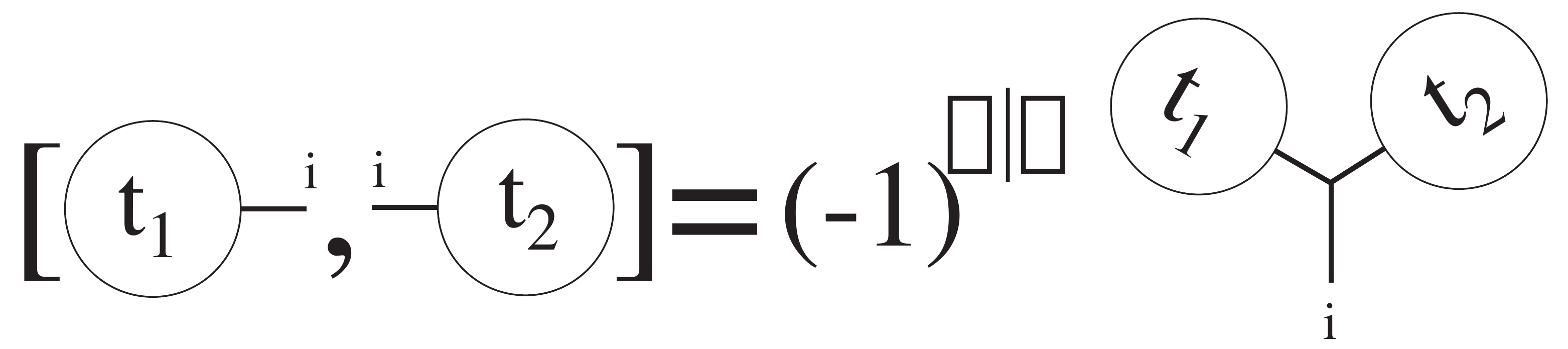}}
\end{definition}
Now we verify that $\mathfrak A_n$ really is a Lie algebra.

\begin{lemma}
$\mathfrak A_n$ is a graded Lie algebra, where the grading is given by the number of trivalent vertices plus one.
\end{lemma}
\begin{proof}
We must show that the bracket satisfies the Antisymmetry and Jacobi identities.
Let $\alpha_i$ be the set of indices involved in $t_i$. 
Antisymmetry is clear, using the identity that $\alpha_1|\alpha_2+\alpha_2|\alpha_1 = |\alpha_1||\alpha_2|$.

Here is the Jacobi identity:
$$[t_1,[t_2,t_3]]=[[t_1,t_2],t_3] +(-1)^{|t_1||t_2|}[t_2,[t_1,t_3]]$$
There are two cases. Either $\alpha_1$, $\alpha_2$ and $\alpha_3$ meet in a single index or 
two of them meet in a single index and the third one meets one of the first two in a different index. If neither of these two cases hold, all three two-fold brackets are zero, and so the equation trivially holds. If they meet in a single index, then each of the three terms represents attaching a tree with two trivalent vertices to the three leaves of $t_1,t_2,t_3$ labeled with the common index $i$, and labeling the left-over leaf $i$. This exactly corresponds to the Jacobi identity with the possible exception of signs.
The sign in front of $[t_1,[t_2,t_3]]$ is $(-1)^{\alpha_2|\alpha_3+\alpha_1|\alpha_2\cup \alpha_3}$,
the sign in front of $[[t_1,t_2],t_3]$ is $(-1)^{\alpha_1|\alpha_2+\alpha_1\cup \alpha_2|\alpha_3}$.
Finally the sign in front of $[t_2,[t_1,t_3]]$ is
$(-1)^{\alpha_1|\alpha_3+\alpha_1\cup \alpha_3|\alpha_2}$.
The first two signs are both equal to $(-1)^{\alpha_1|\alpha_2 +\alpha_2|\alpha_3+\alpha_1|\alpha_3}$. This cancels with some of the third sign, the residue of which is $(-1)^{\alpha_3|\alpha_2+\alpha_2|\alpha_3}=(-1)^{|\alpha_1||\alpha_2|}.$ This gives the correct signs for the standard Jacobi identity. 

The second case is when $\alpha_1,\alpha_2$ and $\alpha_3$ do not meet in a point.
For specificity assume that  $\alpha_1$ and $\alpha_2$ meet in the index $i$ and $\alpha_2$ and $\alpha_3$ meet in the index $j$. Then the term $[t_2,[t_1,t_3]]$ is zero. The other two terms involve splicing $t_2$ to $t_1$ and $t_3$, and so are equal, and even have the same sign, as we argued above.
\end{proof}

\begin{definition}
Let $\bar\A^I_{1,n}$ be defined like $\A^I_{1,n}$ except that we don't divide by 
$\operatorname{STU^2}$ relations.
\end{definition}

Notice that $\bar\A^I_{1,n}\otimes\mathbb Q\subset \mathfrak A_n$ as the top degree term. 

\begin{definition}
We define a map of Lie algebras $\Delta\colon \mathcal B^e_n\to\mathfrak A_n$  which is given on generators by letting $\Delta(x_{ij})$ be the line segment with one end numbered $i$ and the other $j$ with the sign $(-1)^{\{i\}|\{j\}}$, unless $i=j$ when we define $\Delta(x_{ii})=0$.
\end{definition}

Before we show this is well-defined, we give an example.

{\bf Example:} 
\begin{align*}
\Delta([[x_{42},x_{23}],x_{13}])&=-[[4\text{---}2,2\text{---}3],1\text{---}3]\\
&=[\phantom{\texttt{Y}}^4\underset{2}{\texttt{Y}}^3,1\text{---}3]\\
&=\phantom{\texttt{\Large H}}^4_2\texttt{\Large H}^1_3\\
&=\includegraphics[width=.5in]{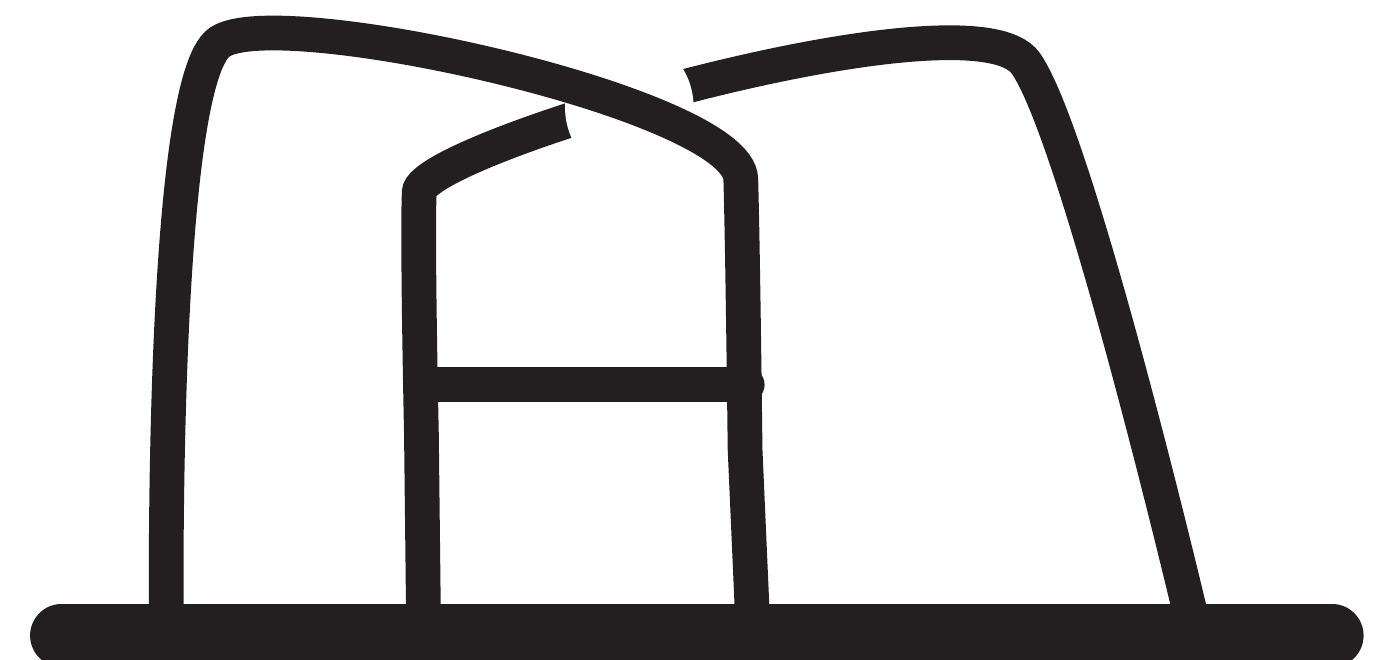}
\end{align*}

\begin{proposition} $\Delta$ is well-defined.
\end{proposition}
\begin{proof}
We must show that $\Delta$ respects the relations from the beginning of Section \ref{overview}.
Clearly $\Delta(x_{ij})=-\Delta(x_{ji})$, since $(-1)^{\{i\}|\{j\}}=-(-1)^{\{j\}|\{i\}}$. Also $\Delta(x_{ii})=0$ by definition. 
Relation $3$ gets sent to zero because $[i\text{---}j,l\text{---}m]\in\mathfrak A_n$ is zero if there are no common indices.
The trees in relation $4$ all get mapped to a trivalent tree with one trivalent vertex and three univalent ones (a ``Y"), with leaves labeled $i,j,l$. The signs work out correctly also, for consider the brackets $[x_{ij},x_{jl}],[x_{li},x_{ij}],[x_{jl},x_{li}]$. 
The vertex orientation of the resulting $Y$'s is the same in all three cases since the cyclic order of indices is the same in all cases. Furthermore, the sign in all three cases is $(-1)^{\{i\}|\{j\}+\{j\}|\{l\}+\{i\}\{l\}}.$ 
\end{proof}

\begin{proposition}
 $\Delta$ induces an isomorphism from $M_{n,n+1}$ to $\bar\A^I_{1,n}\otimes\mathbb Q$.
 \end{proposition}
\begin{proof} 
 Let $\nabla\colon\bar\A^I_{1,n}\to M_{n,n+1}$ be defined as follows. Given a tree $t$, labeled by $1,\ldots,n$, think of $1$ as the root. For each trivalent vertex, there are two index sets, $\alpha$ and $\beta$ representing the numbers labeling the two branches growing away from the root. For each such vertex consider the sign which is $(-1)^{\alpha|\beta}$. Now $\nabla(t)$ is defined to be the iterated commutator formed by replacing each index $k$ by $x_{1k}$, and interpreting the trivalent vertices of $t$ as brackets, multiplied by the product of signs coming from each vertex.
 
One must check that $\nabla$ is well-defined, which means that it respects $\operatorname{IHX}$ (Jacobi identities) and the antisymmetry relations, which is straightforward. Now one checks that $\Delta\circ\nabla=\operatorname{id}_{\bar\A^I_{1,n}}.$  The fact that  $\nabla\circ\Delta=\operatorname{id}_{M_{n,n+1}}$ can be verified most easily if one checks the equality for the generating set of $M_{n,n+1}$ consisting of iterated brackets in the Lie algebra generators $x_{1i}$, $1\leq i\leq n$.
  \end{proof}

\begin{proposition}
The image of $d$ is the subspace of $\operatorname{STU}^2$ relations.
\end{proposition}
\begin{proof}
We may take the domain of $d$, $M_{n,n}$, to be generated by iterated commutators, $c$, in $x_{1i}$ with a single repeated index $i=k$. 

I claim that $M_{n,n}$ is generated by such iterated commutators, with the additional property that $c=[c_1,c_2]$ where $c_1$ and $c_2$ both involve the repeated index $k$. To see this, think of $c$ as a rooted tree in the usual way, and draw a geodesic in the tree connecting the two leaves labeled $x_{1k}$. We want to rewrite the tree so that the root is distance $1$ from this geodesic. Using the Jacobi identity ($\operatorname{IHX}$ relation), one can rewrite a tree as a linear combination of two trees where the two new trees have the geodesic one closer to the root. Continue inductively.

Let $c_1$ be an iterated commutator in the non-repeated generators $x_{1k},x_{1n_1},\ldots, x_{1n_\ell}$, and let $c_2$ be an iterated commutator in the non-repeated generators $x_{1k},x_{1m_1},\ldots,x_{1m_s}$. Let $n_I=\{n_1,\ldots, n_\ell\}$ and $m_J=\{m_1,\ldots, m_s\}$, and assume that $n_I\cup m_J\cup\{1,k\}=\{1,\ldots, n\}$.
Then by the preceding paragraph elements of the form $[c_1,c_2]$ form a generating set for $M_{n,n}$.

Let $c_i[1]$ represent the iterated commutator $c_i$ where each $x_{1j}$ is replaced by
$x_{1,j+1}$, and let $c_i[2]$ represent the iterated commutator where each $x_{1j}$ is replaced by $x_{2,j+1}$. Similarly let $c_i\{k\}$ represent the commutator where each $x_{1j}$ is replaced by $x_{1j}$ if $j<k$, is replaced by $x_{1,j+1}$ if $j>k$ and by
$x_{1k}$ if $j=k$. Also let $c_i\langle k+1\rangle$  represent the commutator where each $x_{1j}$ is replaced by $x_{1j}$ if $j<k$, is replaced by $x_{1,j+1}$ if $j>k$ and by
$x_{1,k+1}$ if $j=k$.

Let $c=[c_1,c_2]$.Now we claim that 
$$d(c)=-\left[c_1[1],c_2[2]]-[c_1[2],c_2[1]\right]+(-1)^k\left[c_1\{k\},c_2\langle k+1\rangle\right]+(-1)^k\left[c_1\langle k+1\rangle,c_2\{k\}\right]$$

To see this,
note that $d(c)=-\tilde\partial^1(c)+(-1)^k\tilde\partial^k(c)$, because no other indices $i$ are repeated, implying the operator $\tilde\partial^i$ is trivial.
$\tilde\partial^1(c)$ has the effect of sending all indices $i>1$ to $i+1$ and summing over changing each $1$ index to either a $1$ or a $2$ in all possible ways where both $1$ and $2$ occur. 
Let us calculate the first term.  Consider a term in  $\bar\partial^1(c)$ where some of the $1$'s have been converted to $2$'s. If there are some $1$'s that remain in $c_1$ and some that remain in $c_2$ then $c_1$ and $c_2$ would have two indices in common, and therefore $c$ would get mapped to zero by $\Delta$. 
So all of the $1$'s need to sit in either $c_1$ or $c_2$ and all of the $2$'s need to sit in the other one. Thus $\tilde\partial^1(c)=[c_1[1],c_2[2]]+[c_1[2],c_2[1]]$.
Similarly $\tilde\partial^k(c)$ converts $k$'s to either $k$'s or $k+1$'s,
which is precisely  $[c_1\{k\},c_2\langle k+1\rangle]+[c_1\langle k+1\rangle,c_2\{k\}]$.

Now we proceed to calculate $\Delta(d(c))$. We claim that $\Delta(c_1[1]),\Delta(c_1[2]),\Delta(c_1\{k\})$ and $\Delta(c_1\langle k+1\rangle)$ are all the same signed tree $t_1$, with leaves labeled by different indices depending on the value of $i$.
To see this, observe that for each, there is an order preserving correspondence between the index sets of the variables $x_{ij}$, and so $\Delta$ assembles them into trees in the same way.
Thus we can write 
\begin{align*}
\Delta(c_1[1])&=t_1(1,k+1,\sigma^1(n_I))\\
\Delta(c_1[2])&=t_1(2,k+1,\sigma^1(n_I))\\
\Delta(c_1\{k\})&=t_1(1,k,\sigma^k(n_I))\\
\Delta(c_1\langle k+1\rangle)&=t_1(1,k+1,\sigma^k(n_I))
\end{align*}
where the indices in parentheses after the tree $t_1$ represent the indices labeling its leaves. 
Similarly, for some signed tree $t_2$,
 \begin{align*}
\Delta(c_2[1])&=t_2(1,k+1,\sigma^1(m_J))\\
\Delta(c_2[2])&=t_2(2,k+1,\sigma^1(m_J))\\
\Delta(c_2\{k\})&=t_2(1,k,\sigma^k(m_J))\\
\Delta(c_2\langle k+1\rangle)&=t_2(1,k+1,\sigma^k(m_J))
\end{align*}

Now
\begin{align*}
\Delta(&-[c_1[1],c_2[2]]-[c_1[2],c_2[1]])=-(-1)^{[2]\cup \sigma^1(n_I)|\{2\}\cup \sigma^1(m_J)}\times\\
&\left([t_1(1,k+1,\sigma^1(n_I)),t_2(2,k+1,\sigma^1(m_J))]-[t_1(2,k+1,\sigma^1(n_I)),t_2(1,k+1,\sigma^1(m_J))]\right)
\end{align*}
and
\begin{align*}
\Delta((-1)^k[&c_1\{k\},c_2\langle k+1\rangle]+(-1)^k[c_1\langle k+1\rangle,c_2\{k\}])=
(-1)^{k+\left(\langle k+1\rangle\cup \sigma^k(n_I)|\{k+1\}\cup \sigma^k(m_J)\right)}\times\\
&\left([t_1(1,k,\sigma^k(n_I)),t_2(1,k+1,\sigma^k(m_J))]-[t_1(1,k+1,\sigma^k(n_I)),t_2(1,k,\sigma^k(m_J))]
\right)
\end{align*}

In the first two terms we are splicing the trees together along the leaves corresponding to the second slot, and labeling the first-slot leaves with $1$ and $2$ in both orders. 
In the second two terms we are splicing $t_1$ and $t_2$ together along the leaves corresponding to the first slot in the parentheses, and the leaves corresponding to the second slot are labeled with $k$ and $k+1$ in both orders. Thus the sum of the four terms is an $\operatorname{STU^2}$ relator, up to signs. To see that the signs work out correctly, note the following identities 
\begin{align*}
[2]\cup \sigma^1(n_I)|\{2\}\cup \sigma^1(m_J)&=|n_I|+ n_I|m_J\\
\langle k+1\rangle\cup \sigma^k(n_I)|\{k+1\}\cup \sigma^k(m_J)&=\langle k+1\rangle|m_J+n_I|\langle k+1\rangle+n_I|m_J\\
\langle k+1\rangle|n_I+n_I|\langle k+1\rangle&=|n_I|\\
\langle k+1\rangle|n_I+\langle k+1\rangle|m_J&=k-2
\end{align*}
The first three are fairly straightforward. To see the last identity, note that the left hand side is equal to the number of indices from $n_I\cup m_J$ which are less than $k$. Since $n_I\cup m_J$ hits everything except for $1$, we count $k-2$. 
With these identities, we see that the signs in front of each pair of terms above is opposite, exactly as needed for an $\operatorname{STU^2}$ relator.

In this way we realize all $\operatorname{STU}^2$ relations where one of the ``splitting sites" is all the way to the left on the line segment. However any $\operatorname{STU}^2$ relation can be written as a difference of two $\operatorname{STU}^2$-relations of this type. 
\end{proof}

The proof of Theorem~\ref{main} is now complete, since $\Delta$ induces an isomorphism
$$\Delta\colon M_{n,n+1}/im(d)\cong A^I_{1,n}\otimes\mathbb Q.$$

\section{Appendix}

In this section we verify that $\Psi'\colon \A^I_{n,n}\to\A^I_{n-2,n}$ vanishes on an ``initial" $4T$ relator.

To facilitate our calculations, we need a lemma. Recall the heavy grey line convention from Figure~\ref{psidef}.

\begin{lemma}\label{idunno}

\noindent\begin{enumerate}
\item The following identity is a consequence of $\operatorname{STU}^2$ relations.
\begin{center}
\includegraphics[width=.7\linewidth]{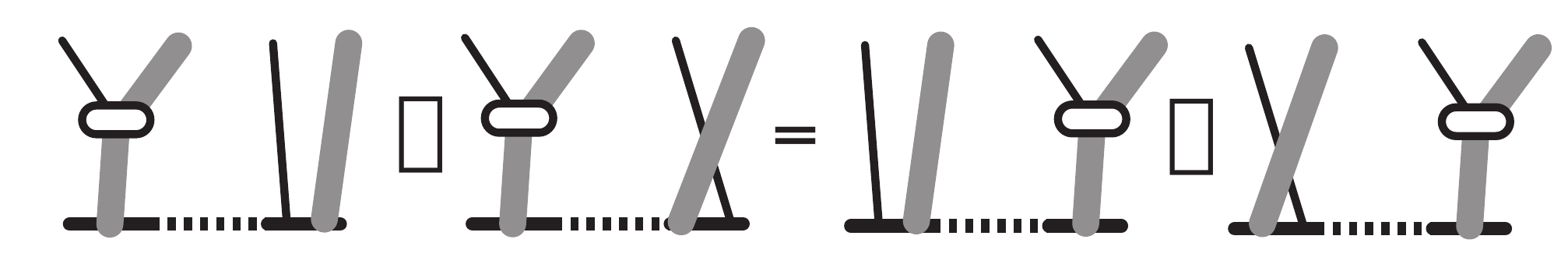}
\end{center}
\item The following identity is a consequence of $\operatorname{IHX}$ relations.
\begin{center}
\includegraphics[width=.5\linewidth]{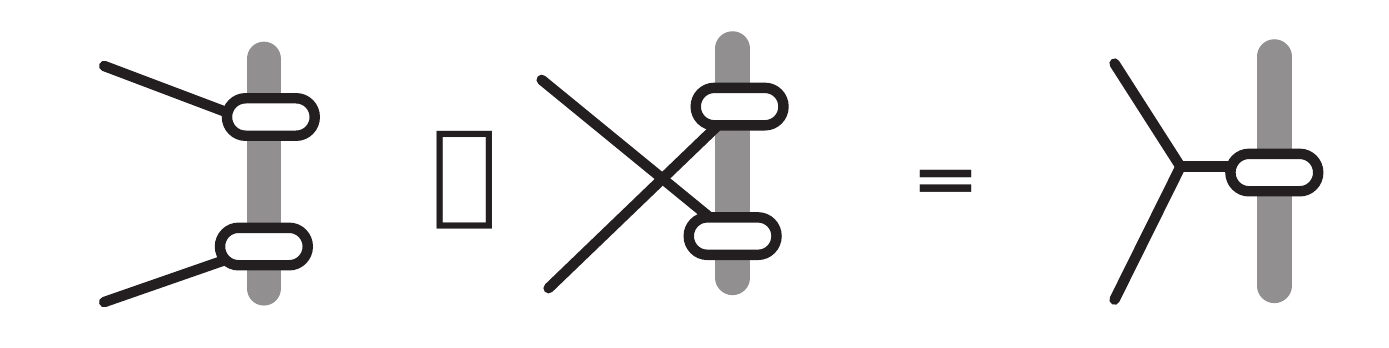}
\end{center}
\end{enumerate}
\end{lemma}
\begin{proof}
The first part is straightforward. Each term represents a linear combination of distributing the end of the black edge touching the oval to all of the edges coming through the grey area. Expanding out all four terms above, they naturally group into $\operatorname{STU}^2$ relations.

The proof of the second part is not much harder, and is encapsulated in Figure~\ref{ihxproof}. The heavy grey line is replaced by the pieces of edges that run through it.
\begin{figure}[ht!]
\begin{center}
\includegraphics[width=.7\linewidth]{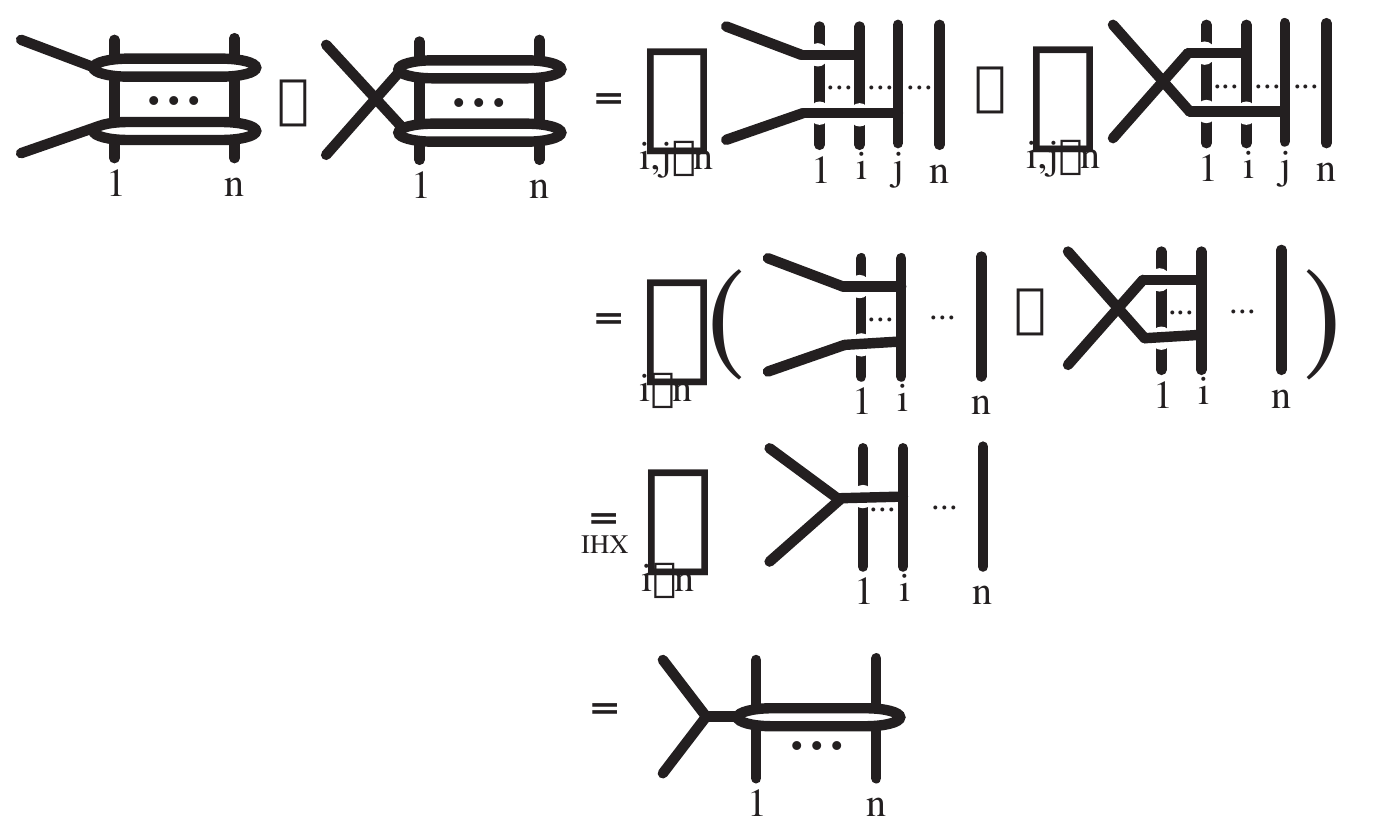}
\caption{The proof of Lemma~\ref{idunno} (2).}\label{ihxproof}
\end{center}
\end{figure}
\end{proof}

\begin{figure}[ht!]
\begin{center}
\includegraphics[width=.85\linewidth]{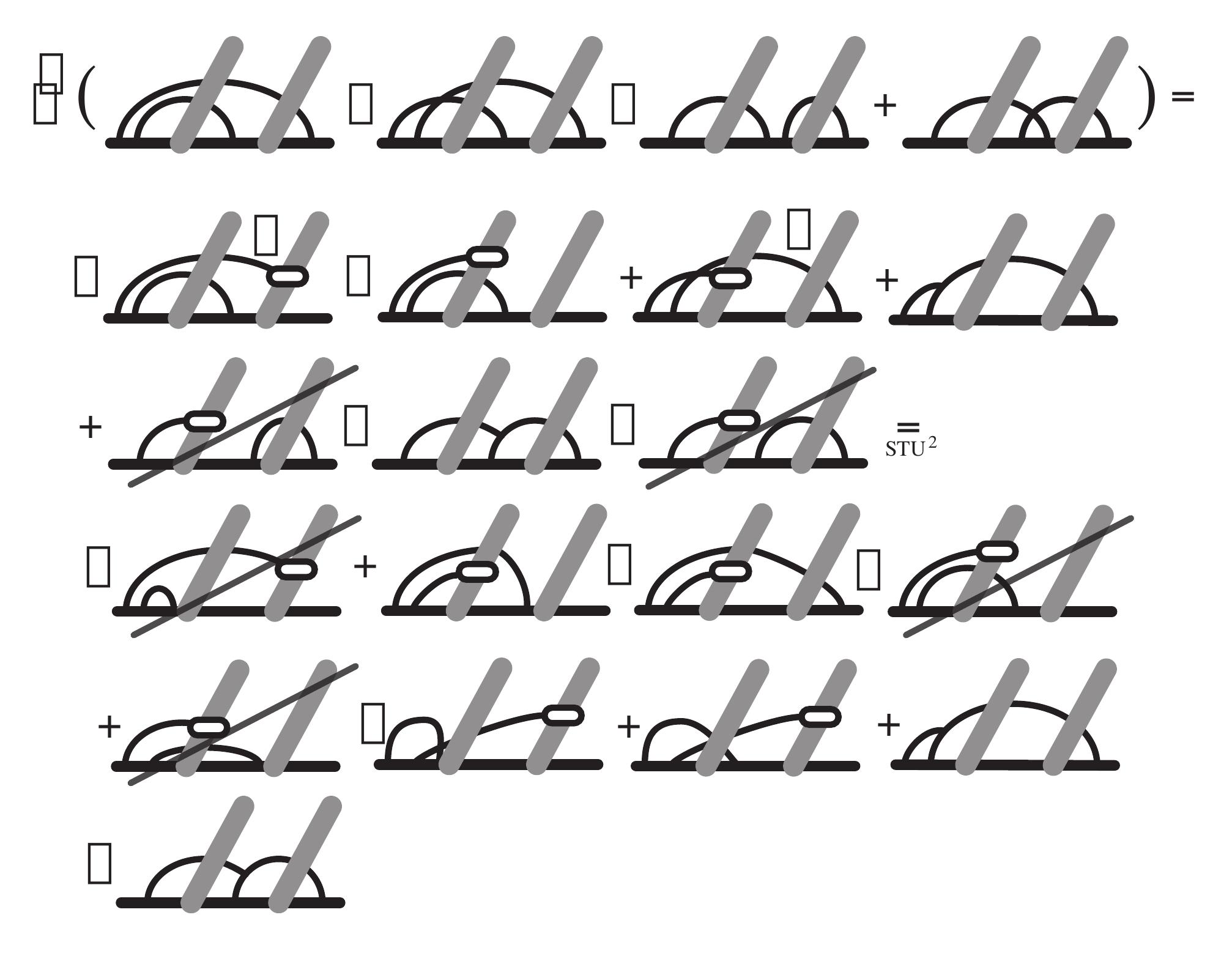}

\includegraphics[width=.85\linewidth]{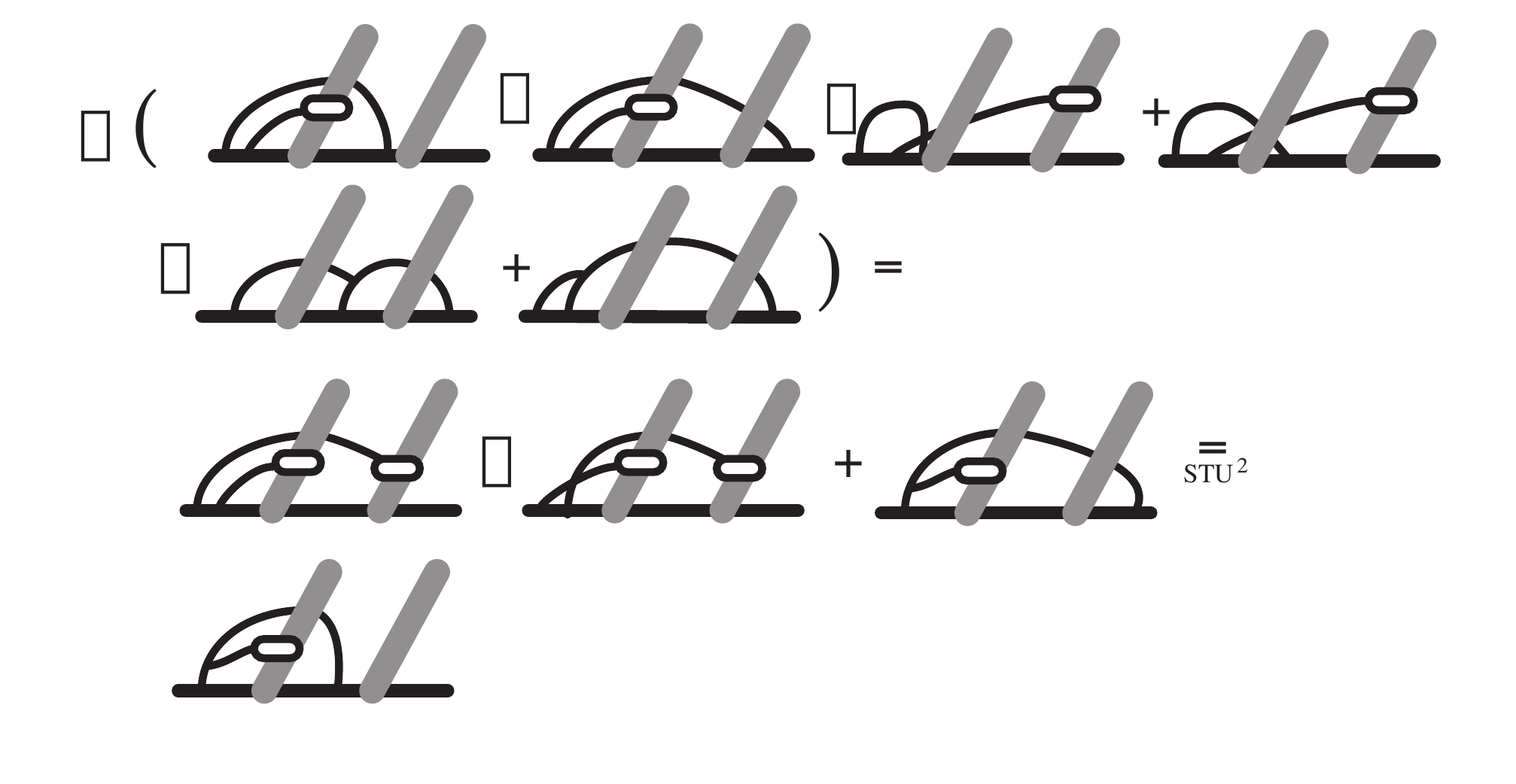}

\caption{Top: Applying $\tilde{\Psi}$ to an ``initial" $4T$ relator. The second equality
is an equivalence modulo $\operatorname{STU}^2$ relations involving the terms with asterisks.
Bottom: Applying $\Psi$ again. The result is $\Psi'$ of an initial $4T$ relator.
}\label{psi1}
\end{center}
\end{figure}

\begin{lemma}\label{vanish}
$\Psi'$ vanishes on $4T$ relations involving a chord which attaches farthest to the left.
\end{lemma}
\begin{proof}
The $4T$ relation is depicted in Figure~\ref{blash}a. 
Let us first consider a relation of the form coming from the first equality in Figure~\ref{blash}a.  
At the top of Figure~\ref{psi1} such a $4T$ relator, involving the left-most chord, is depicted, with heavy grey lines indicating places in the diagram where many edges from a tree may hit, some of which may have multiple endpoints inside the grey region. 
Then $\tilde \Psi$ is applied. The result is reduced modulo $\operatorname{STU}^2$ and $\operatorname{SEP}$ relations, and then $\Psi$ is applied again. One can use Claim~\ref{clm}
to calculate this modulo $\operatorname{STU}^2$ relations, although one can also verify that equality holds on the nose. The result, $\Psi'$ of the $4T$ relator, is equivalent modulo $\operatorname{STU}^2$ relations to the picture at the bottom. It now remains to show that this linear combination of diagrams is trivial modulo $\operatorname{STU}^2$ and $\operatorname{IHX}$ relations.The two equations in Figure~\ref{adhoc} show that this is true,
and these equations follow from Lemma \ref{idunno}.

\begin{figure}
\begin{center}
\includegraphics[width=.8\linewidth]{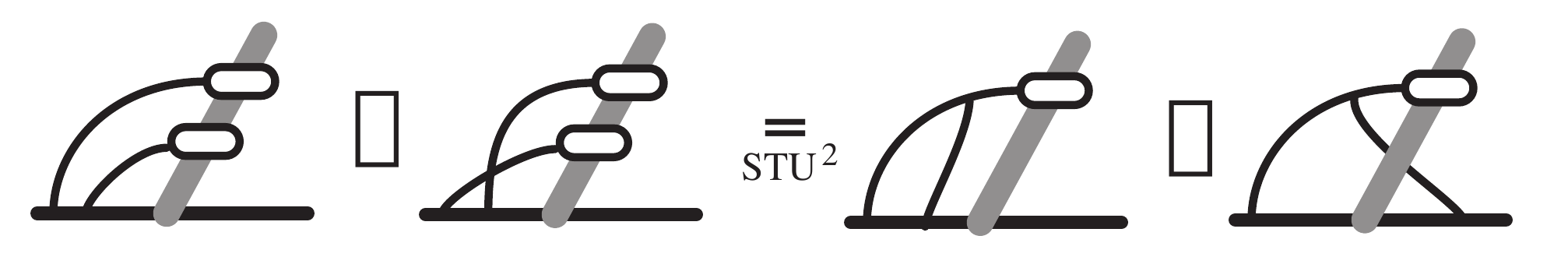}\\
\includegraphics[width=.8\linewidth]{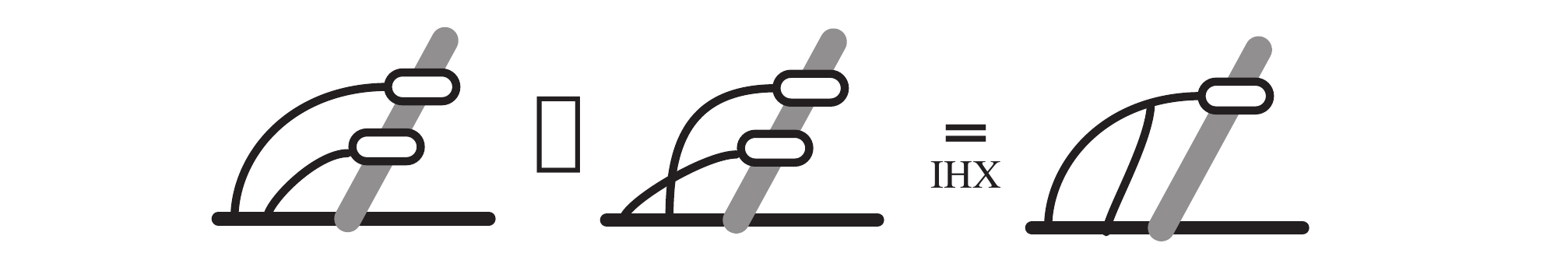}
\caption{From the proof of Lemma~\ref{vanish}.}\label{adhoc}
\end{center}
\end{figure}

Finally, we need to consider a relation from Figure~\ref{blash}a equating the first and third pairs of terms. A straightforward calculation shows that $\tilde{\Psi}$ already sends this to zero. The calculation is simpler, in this case, because the ``same" chord is farthest to the left in all four terms.
\end{proof}

\end{document}